\documentclass{amsart}
\usepackage{amsmath}
\usepackage{amsfonts}
\usepackage{hyperref}

\newtheorem{theorem}{Theorem}
\theoremstyle{plain}

\newtheorem{claim}{Claim}

\newtheorem{corollary}{Corollary}

\newtheorem{definition}{Definition}
\newtheorem{example}{Example}

\newtheorem{lemma}{Lemma}

\newtheorem{problem}{Problem}
\newtheorem{proposition}{Proposition}
\newtheorem{remark}{Remark}

\numberwithin{equation}{section}
\input{tcilatex}

\newcommand{\defin}[1]{\textbf{\emph{#1}}}

\newcommand{\R}{\mathbb{R}}

\newcommand{\diam}{\mathrm{diam}}
\newcommand{\norm}[2]{\left\Vert{#1}\right\Vert_{#2}}
\newcommand{\dom}{\mathrm{dom}}
\newcommand{\supp}{\mathrm{Supp}}

\newcommand{\ceil}[1]{\lceil{#1}\rceil}
\newcommand{\integral}[2]{\int\!{#1}\,\mathrm{d}{#2}}
\newcommand{\abs}[1]{\left|{#1}\right|}

\begin{document}
\title[Constructing pseudo-random points]{A constructive Borel-Cantelli Lemma.
Constructing orbits with required statistical properties.\footnote{partly
  supported by ANR Grant  05 2452 260 ox}}
\date{30 juin 2008}
\author{Stefano Galatolo}
\address{Dipartimento di Matematica Applicata, Universita di Pisa}
\email{s.galatolo@docenti.ing.unipi.it}
\author{Mathieu Hoyrup}
\address{LIENS, Ecole Normale Sup\'erieure, Paris}
\email{hoyrup@di.ens.fr}
\author{Crist\'obal Rojas}
\address{LIENS, Ecole Normale Sup\'erieure and CREA, Ecole Polytechnique,
Paris}
\email{rojas@di.ens.fr}

\begin{abstract}
In the general context of computable metric spaces and computable measures
we prove a kind of constructive Borel-Cantelli lemma: given a sequence
(constructive in some way) of sets $A_{i}$ with effectively summable measures,
there are computable points which are not contained in infinitely many $A_{i}$.

As a consequence of this we obtain the existence of computable points which
follow the \emph{typical statistical behavior} of a dynamical system (they
satisfy the Birkhoff theorem) for a large class of systems, having
computable invariant measure and a certain ``logarithmic'' speed of convergence of Birkhoff averages over Lipshitz observables. This is
applied to uniformly hyperbolic systems, piecewise expanding maps, systems
on the interval with an indifferent fixed point and it directly implies the
existence of computable numbers which are normal with respect to any base.
\end{abstract}

\maketitle
\tableofcontents


\section{Introduction}

Many results in mathematics ensure the existence of points satisfying a
given property $P$ by estimating the measure of $P$ and proving that it is
positive. In general this approach is not constructive and does not give an
effective way to construct points satisfying the given property.

A key lemma in this kind of techniques is the well-known Borel-Cantelli one:

\noindent \textbf{Borel-Cantelli Lemma. }\emph{Let }$\{A_{n}\}$\emph{\ be a
sequence of subsets in a probability space }$(X,\mu )$\emph{. If }$\sum \mu
(A_{n})<\infty $\emph{, then }$\mu (\lim \sup A_{n})=0,$\emph{\ that is, the
set of points which are contained in infinitely many }$A_{n}$\emph{\ has
zero measure}.

Under these conditions, $X-\lim \sup A_{n}$ is a full measure set and hence
it contains ``many'' points of $X$. In this paper we give a general method
to construct points in this set. This method will be applied to some
nontrivial problems, as constructing numbers which are normal in every base
and typical trajectories of dynamical systems.

To face this problem we will put ourself in the framework of computable
metric spaces. Let us introduce and motivate this concept. It is well known
that the state of a physical system can be \emph{known} only up to some
finite precision (because of measuring errors, thermal shaking, quantum
phenomena, long range interactions etc...). From a mathematical point of
view this knowledge is represented by a ball with positive radius in the
metric space of all possible configurations of the system.

In practice, the knowledge of the state of the system up to some finite
precision can be described by a sentence like \textquotedblleft the position
of the point in the phase space at time $3$ is $x(3)=0.322\pm 0.001$%
\textquotedblright . What is important here is that it admits a \emph{finite
description} (a finite string of characters).

This finite string of characters, can then be elaborated to estimate, for
example the position or the distance of the system's status at time 3 with
respect to other points of the space.

This kind of identification 
\begin{equation*}
Strings\leftrightarrow[Points, Geometrical\ objects]
\end{equation*}
if often implicit, and considered to be obvious but it underlies the concept
of Computable Metric Space.

A Computable Metric Space is a metric space where a dense countable set
(which will be called the set of ideal points) is identified with a set of
finite strings, in a way that the distance between points in this set can be
computed up to any given approximation by an algorithm having the
corresponding strings as an input (see section \ref{CMS} for precise
definitions).

For example in $\mathbb{R}$ the set $\mathbb{Q}$ can be identified with the
strings ``$p$ over $q$'' in a way that the distance between rationals can be
obviously calculated by an algorithm having the strings as input. We remark
that if $\mathbb{R}$ is considered as a computable metric space, then beyond 
$\mathbb{Q}$ there are many other points which admit finite descriptions,
for example $\pi $ or $\sqrt{2}$ are not rationals but they can be
approximated at any given precision by an algorithm, hence in some sense
this points too can be identified to finite strings: $\pi $ for example can
be identified with the finite program which approximates it by rationals at
any given precision. This set of points is called the set of \emph{computable} real numbers (they were introduced by Turing in \cite{Tur36}).
The concept of computable point can be easily generalized to any computable
metric space. Coming back to our main question, now the problem we face is
the following: Given some property $P$ about points of $X$ (or equivalently
a subset of $X$), can this property be observed with a computer? That is,
does there exist computable points satisfying this property?

For instance, given a (non atomic) probability measure $\mu $, let $P$ be a
subset of $X$ of probability one: a point chosen \textquotedblleft at
random\textquotedblright\ will almost surely belong to $P$. But, as the set
of computable points have null measure (is a countable set) the full measure
of $P$ induces \emph{a priori} nothing upon its computable part (i.e. the
set of computable points belonging to $P$).

We will give some results which give a positive answer to this question when 
$P$ is constructed by a Borel-Cantelli technique. Let us illustrate this
(for a precise statement see theorem \ref{effective_BC_theorem}):\newline

\noindent \textbf{Theorem A. }\emph{Let us consider a sequence of closed
sets $(A_{n})_{n\in \mathbb{N}}$ (with some effectivity condition, see
definition \ref{reop}) such that $\sum \mu (A_{n})<\infty $ in an effective
way (see Def. \ref{ebc}).}

\emph{If the
measure $\mu$ is computable (Def. \ref{compmeas} ) then there are computable
points outside $\limsup A_n$, that is lying in $A_n$'s only finitely many
times.}\newline

\noindent \textbf{Computable absolutely normal numbers}. As an example, a
classical question where this kind of tool can be naturally applied is the
normality: given a fixed enumeration base $b$ of real numbers it it quite
easy to prove that the set of $b$-normal numbers (the numbers where all the
digits $\{0,...,b-1\}$ appear with the same frequency) has Lebesgue-measure
one. Can we find computable normal numbers? The construction proposed by
Champernowne \cite{Cha33} happens to be algorithmic, so it gives a positive
answer to the question.

A natural and much more difficult problem is to construct numbers which are
normal in \emph{every} base (see sec. \ref{normal} for some historical
comments on the problem). In section \ref{normal} the existence of
computable absolutely normal numbers will be obtained as a quite simple
corollary of Theorem A.

\noindent \textbf{Computable points having typical statistical behavior}.
The above result on normal numbers is a particular case of the construction
of computable points which follows the typical statistical behavior of a
dynamical system. We will need the notion of computable dynamical systems,
let us introduce it.

The notion of algorithm and computable function can be extended to functions
between computable metric spaces (Def. \ref{comp_func}). This allows to
consider computable dynamical systems over metric spaces (systems whose
dynamics is generated by the iteration of a computable function), and
computable observables. With these definitions, all systems which can be
effectively described (and used in simulations) are computable.

Computable points (as described above) are a very small invariant set,
compared to the whole space. By this reason, a computable point rarely can
be expected to behave as a typical point of the space and give rise to a
typical statistical behavior of the dynamics. Here, ``typical'' behavior
means a behavior which is attained for a full measure set of initial
conditions. Nevertheless computable points are the only points we can use
when we perform a simulation or some explicit computation on a computer.

A number of theoretical questions arise naturally from all these facts. Due
to the importance of the general forecasting-simulation problem these
questions also have a practical importance.

\begin{problem}
Since simulations can only start with computable initial conditions, given
some typical statistical behavior of a dynamical system, is there some
computable initial condition realizing this behavior? how to choose such
points?
\end{problem}

Such points could be called \emph{pseudorandom} points. Meaningful
simulations, showing typical behaviors of the dynamics can be performed if
computable, pseudorandom initial conditions exist\footnote{It is widely believed that computer simulations produce correct ergodic
behaviour. The evidence is mostly heuristic. Most arguments are based on the
various ``shadowing'' results(see e.g. \cite{KH} chapter 18).
In this kind of approach (different from our), it is possible to prove that
in a suitable system, any ``pseudo''
-trajectory, as the ones which are obtained in simulations with some
computation error is near to a real trajectory of the system.
\par
So we know that what we see in a simulation is near to some real trajectory
(even if we do not know if the trajectory is typical in some sense). The
main limit of this approach is however that shadowing results hold only in
particular systems, having some uniform hyperbolicity, while many physically
interesting systems are not like this.
\par
We recall that in our approach we consider real trajectories instead of
``pseudo'' ones and we ask if there is some
computable point which behaves as a typical point of the space.}. A somewhat similar
problem has already been investigated in \cite{KruTro94} in the setting of
symbolic dynamics. They consider recursive discretisations of the system
(that is a subset of computable points) and look for conditions to ensure
that a \emph{finite observer} is unable to distinguish the motion on the
recursive discretisation from the original system.

In our framework, a first topological result is the following: if the system
is computable and has at least a dense orbit, then there is a computable
point having a dense orbit (see Thm. \ref{computable_dense}).

From the statistical point of view we can use the above Theorem A to prove
the following second main result which we summarize informally below (see
Thm.s \ref{typical_f} and \ref{mutypicalteo} for precise statements).\newline

\noindent \textbf{Theorem B.} \emph{If \ }$(X,\mu,T)$\emph{\ is a computable
dynamical system and}

\begin{enumerate}
\item $\mu $\emph{\ is a computable invariant ergodic measure.}

\item \emph{The system }$(X,T,\mu )$\emph{\ is $log^2$  ergodic (see definition \ref{logarithmically ergodic} ) for observables in some suitable functions space.}
\end{enumerate}

\emph{Then there exist computable points $x$ for which it holds:} 
\begin{equation}  \label{typic}
\underset{n\rightarrow \infty }{\lim }\frac{f(x)+f(T(x))+\ldots
+f(T^{n-1}(x))}{n}=\integral{f}{\mu}
\end{equation}%
\emph{for any continuous function $f:X\to \mathbb{R}$ with compact support.}
\newline

The above theorem states that in such systems there are computable points
whose time average equals the space average for any such observable on $X$,
hence providing a set of computable points which from the statistical point
of view behave as the typical points of $(X,\mu)$ in the Birkhoff pointwise
ergodic theorem.

We remark that the approach taken in \cite{KruTro94} is quite different,
in the sense that they give sufficient conditions (in terms of Kolmogorov
complexity) for a subset of computable points (a recursive discretisation)
which ensure that this set satisfies a kind of finite ergodic theorem (a
much weaker property than \ref{typic}) but give no method to construct such
computable points (because these conditions cannot be verified in a
constructive way).

To apply theorem B to concrete systems the main difficulty is to verify the
points 1) and 2). In section \ref{examples} we show that these are verified
for the SRB invariant measure (the natural invariant measure to be
considered in this cases) in some classes of interesting systems as
uniformly hyperbolic systems, piecewise expanding maps and interval maps
with an indifferent fixed point.


\section{Computability}

The starting point of recursion theory was to give a mathematical definition
making precise the intuitive notions of \emph{algorithmic} or \emph{%
effective procedure} on symbolic objects. Every mathematician has a more or
less clear intuition of what can be computed by algorithms: the
multiplication of natural numbers, the formal derivation of polynomials are
simple examples.

Several very different formalizations have been independently proposed (by
Church, Kleene, Turing, Post, Markov...) in the 30's, and have proved to be
equivalent: they compute the same functions from $\mathbb{N}$ to $\mathbb{N}$%
. This class of functions is now called the class of \emph{recursive
functions}. As an algorithm is allowed to run forever on an input, these
functions may be \emph{partial}, i.e not defined everywhere. The \emph{domain%
} of a recursive function is the set of inputs on which the algorithm
eventually halts. A recursive function whose domain is $\mathbb{N}$ is said
to be \emph{total}.

We now recall an important concept from recursion theory. A set $E\subseteq 
\mathbb{N}$ is called \textbf{\emph{recursively enumerable (r.e.)}} if there
is a (partial or total) recursive function $\varphi :\mathbb{N}\rightarrow 
\mathbb{N}$ enumerating $E$, that is $E=\{\varphi (n):n\in \mathbb{N}\}$. If 
$E\neq\emptyset$, $\phi$ can be effectively converted into a total recursive
function $\psi$ which enumerates the same set $E$. We recall a useful
characterization of r.e. sets: a set $E\subseteq \mathbb{N}$ is said to be 
\textbf{\emph{semi-decidable}} if there is a recursive function $\varphi :%
\mathbb{N}\rightarrow \mathbb{N}$ whose domain is $E$, that is $\varphi (n)$
halts if and only if $n\in E$. A set is r.e. if and only if it is
semi-decidable, and the corresponding recursive functions can be effectively
converted one another. We will freely use this equivalence, using in each
particular situation the most adapted characterization.

\subsection{Algorithms and uniform algorithms}

Strictly speaking, recursive functions only work on natural numbers, but
this can be extended to the objects (thought as ``finite'' objects) of any
countable set, once a numbering of its elements has been chosen. We will use
the word \emph{algorithm} instead of \emph{recursive function} when the
inputs or outputs are interpreted as finite objects. The operative power of
algorithms on the objects of such a numbered set obviously depends on what
can be effectively recovered from their numbers.

More precisely, let $X$ and $Y$ be such numbered sets such that the
numbering of $X$ is injective (it is then a bijection between $\mathbb{N}$
and $X$). Then any \emph{recursive function} $\varphi:\mathbb{N}\to\mathbb{N}
$ induces an \emph{algorithm} $\mathcal{A}:X\to Y$. The particular case $X=%
\mathbb{N}$ will be much used.

For instance, the set $\mathbb{Q}$ of rational numbers can be injectively
numbered $\mathbb{Q}=\{q_0,q_1,\ldots\}$ in an \emph{effective} way: the
number $i$ of a rational $a/b$ can be computed from $a$ and $b$, and vice
versa. We fix such a numbering: from now and beyond $q_i$ will designate the
rational number which has number $i$.

Now, let us consider a computability notion in the real number set, here for
a number to be computable means that there is an algorithm which can
approximate the number up to any error. We remark that this notion was
already introduced by Turing in \cite{Tur36}.

Let $x$ be a real number and define $\mathbb{Q}^{<}(x):=\{i\in \mathbb{N}
:q_i<x\}$.

\begin{definition}
We say that:

\begin{itemize}
\item $x$ is \textbf{\emph{lower semi-computable}} if the set $\mathbb{Q}%
^{<}(x)$ is r.e.

\item $x$ is \textbf{\emph{upper semi-computable}} if the set $\mathbb{Q}%
^{<}(-x)$ is r.e.

\item $x$ is \textbf{\emph{computable}} if it is lower and upper
semi-computable.
\end{itemize}
\end{definition}

Equivalently, a real number is computable if and only if there exists an
algorithmic enumeration of a sequence of rational numbers converging
exponentially fast to $x$. That is:

\begin{proposition}
A real number is \textbf{\emph{computable}} if there is an algorithm $%
\mathcal{A}:\mathbb{N}\to \mathbb{Q}$ such that $|\mathcal{A}(n)-x|\leq
2^{-n}$ for all $n$.
\end{proposition}

\noindent \textbf{Uniformity.}

Algorithms can be used to define computability notions on many classes of
mathematical objects. The definition of computability notions will be
particular to a class of objects, but they will always follow the following
scheme:

\begin{center}
An object $O$ is \textbf{\emph{computable}} if there is an \\[0pt]
algorithm $\mathcal{A}:\to Y$ which computes $O$ in some way.
\end{center}

Each computability notion comes with a uniform version. Let $(O_i)_{i\in%
\mathbb{N}}$ be a sequence of computable objects:

\begin{center}
$O_i$ is computable \textbf{\emph{uniformly in $\boldsymbol{i}$}} if there
is an algorithm \\[0pt]
$\mathcal{A}:\mathbb{N}\times X \to Y$ such that for all $i$, $\mathcal{A}_i=%
\mathcal{A}(i,.):X\to Y$ computes $O_i$.
\end{center}

For instance, the elements of a sequence of real numbers $(x_i)_{i\in\mathbb{%
N}}$ are uniformly computable if there is a algorithm $\mathcal{A}:\mathbb{N}%
\times\mathbb{N}\to\mathbb{Q}$ such that $|\mathcal{A}(i,n)-x_i|\leq 2^{-n}$
for all $i,n$.

In each particular case, the computability notion may take a particular
name: computable, constructive, effective, r.e., etc... so the term
``computable'' used above shall be replaced. 

\subsection{Computable metric spaces\label{CMS}}

A computable metric space is a metric space with an additional structure
allowing to interpret input and output of algorithms as points of the metric
space (for an introduction to other approaches to this concept see \cite%
{HoyRoj07}). This is done in the following way: there is a dense subset
(called ideal points) such that each point of the set is identified with a
natural number. The choice of this set is compatible with the metric, in the
sense that the distance between two such points is computable up to any
precision by an algorithm getting the names of the points as input. Using
this simple assumptions many constructions on metric spaces can be
implemented by algorithms.

\begin{definition}
A \textbf{\emph{computable metric space}} (CMS) is a triple $\mathcal{X}%
=(X,d,S)$, where

\begin{enumerate}
\item[(i)] $(X,d)$ is a separable metric space.

\item[(ii)] $S=\{s_i\}_{i\in\mathbb{N}}$ is a countable set of elements from 
$X$ (called \textbf{\emph{ideal points}}) which is dense in $(X,d)$.

\item[(iii)] The distances between ideal points $d(s_{i},s_{j})$ are all
computable, uniformly in $i,j$ (there is an algorithm that gets the names of
two points and an allowed error as an input and outputs the distance between
two points up to the given approximation).
\end{enumerate}
\end{definition}

$S$ is a numbered set, and the information that can be recovered from the
numbers of ideal points is their mutual distances. Without loss of
generality, we will suppose the numbering of $S$ to be injective: it can
always be made injective in an effective way.

We say that in a metric space $(X,d)$, a sequence of points $(x_n)_{n\in%
\mathbb{N}}$ converges \textbf{\emph{fast}} to a point $x$ if $d(x_n,x)\leq
2^{-n}$ for all $n$.

\begin{definition}
A point $x\in X$ is said to be \textbf{\emph{computable}} if there is an
algorithm $\mathcal{A}:\mathbb{N}\to S$ such that $(\mathcal{A}(n))_{n\in%
\mathbb{N}}$ converges fast to $x$.
\end{definition}

We define the set of \textbf{\emph{ideal balls}} to be $\mathcal{B}%
:=\{B(s_{i},q_{j}):s_{i}\in S,q_{j}\in \mathbb{Q}_{>0}\}$. We fix a
numbering $\mathcal{B}=\{B_0,B_1,\ldots\}$ which makes the number of a ball
effectively computable from its center and radius and vice versa (this
numbering may not be injective). $\mathcal{B}$ is a countable basis of the
topology.

\begin{definition}[Constructive open sets]
\label{reop} We say that an open set $U$ is \textbf{\emph{constructive}} if there is an algorithm $\mathcal{A}:\mathbb{N}\to \mathcal{B}$ such that 
$U=\bigcup_n \mathcal{A}(n)$.
\end{definition}

Observe that an algorithm which diverges on each input $n$ enumerates the
empty set, which is then a constructive open set. Sequences of uniformly constructive open sets are
naturally defined.

\begin{example}
We give some example of constructive open sets:

\begin{itemize}
\item The whole space $X$ is constructive open.

\item Every finite union or intersection of ideal balls $\{B_{n_1},%
\ldots,B_{n_k}\}$ is a constructive open set, uniformly in ${\langle n_1,\ldots,n_k
\rangle}$.

\item If $(U_i)_{i\in\mathbb{N}}$ is a sequence of uniformly constructive open sets,
then $\bigcup_i U_i$ is a constructive open set.
\end{itemize}
\end{example}

\begin{remark}
If $U$ is constructively open, belonging to $U$ for an ideal point is semi-decidable:
there is an algorithm $\mathcal{A}:S\to \mathbb{N}$ which halts only on
ideal points belonging to $U$. Equivalently, the set of ideal points lying
in $U$ is r.e. (as a subset of $S$): there is an algorithm $\mathcal{A}:%
\mathbb{N}\to S$ enumerating $S\cap U$. Hence $(U,S\cap U,d)$ has a natural
structure of computable metric space.
\end{remark}

\begin{definition}[Constructive $G_\delta$-set]
A \defin{constructive $\boldsymbol{G_\delta}$-set} is an intersection of a sequence of uniformly constructive open sets.
\end{definition}

Obviously, an intersection of uniformly constructive $G_\delta$-sets is
also a constructive $G_\delta$-set.

Let $(X,S_{X}=\{s_{1}^{X},s_{2}^{X},...\},d_{X})$ and
$(Y,S_{Y}=\{s_{1}^{Y},s_{2}^{Y},...\},d_{Y})$ be computable metric
spaces. Let also $B_{i}^{X}$ and $B_{i}^{Y}$ be enumerations of the ideal
balls in $X$ and $Y$. A computable function $X\rightarrow Y$ is a
function whose behavior can be computed by an algorithm up to any
precision. For this it is sufficient that the preimage of each ideal ball
is calculated with any precision.

\begin{definition}[Computable Functions]\label{comp_func}
A function $T:X\rightarrow Y$ is \defin{computable} if $T^{-1}(B_{i}^{Y})$
is a constructive open set, uniformly in $i$. That is, there is an
algorithm $\mathcal{A}:\mathbb{N}\times \mathbb{N}\rightarrow
\mathcal{B}^{X} $ such that $T^{-1}(B_{i}^{Y}
)=\bigcup_{n}\mathcal{A}(i,n)$ for all $i$.

A function $T:X\rightarrow Y$ is \defin{computable on $D\subseteq X$} if there are uniformly constructive open sets $U_i$ such that $T^{-1}(B_{i}^{Y})\cap D=U_i\cap D$.
\end{definition}

\begin{remark}
\label{tracci}We remark that if $T$ is computable then all 
$T(s_{i}^{X})$ are computable uniformly in $i$: there is an algorithm $\mathcal{A}:\mathbb{N}\times \mathbb{N}\to S^Y$ such that $(\mathcal{A}(i,n))_{n\in\mathbb{N}}$ converges fast to $T(s^X_i)$ for all $i$.

The algorithm just semi-decides for each ideal ball in $Y$ if $s_{i}^{X}$ is
contained in its preimage. The process will stop for each ideal ball that
contains $T(s_{i}^{X})$, which allows to extract a sequence of ideal points
of $Y$ which converges fast to $T(s^X_i)$.
\end{remark}

We remark that there are many more or less equivalent characterizations of
computable functions between CMS, see \cite{HoyRoj07} for more details.

The following (\cite{Gal00}, Def. 21 ) is a criteria to check computability
of a large class of uniformly continuous functions.

\begin{proposition}
\label{unifcont}If $T$ satisfies the following:

\begin{itemize}
\item all $T(s_{i}^{X})$ are computable points, uniformly in $i$,

\item $T$ is recursively uniformly continuous: there is an algorithm $%
\mathcal{A}:\mathbb{Q}_{>0}\rightarrow \mathbb{Q}_{>0}$ such that for all $%
\epsilon \in \mathbb{Q}_{>0}$, $d(x,x^{\prime })<\mathcal{A}(\epsilon
)\Rightarrow d(T(x),T(x^{\prime }))<\epsilon $,
\end{itemize}
then $T$ is computable.
\end{proposition}

\begin{proof}
Let $E=\{(i,j)\in\mathbb{N}^2: d(T(s^X_i),s)+q_j<r\}$: this is a r.e. subset
of $\mathbb{N}$ (uniformly in $s,r$) by the first condition. Then one can
show that $T^{-1}(B(s,r))=\bigcup_{(i,j)\in E} B(s_i,A(q_j))$.
\end{proof}




\subsection{Computable measures\label{seccompmu}}

When $X$ is a computable metric space, the space of probability measures
over $X$, denoted by $\mathcal{M}(X)$, can be endowed with a structure of
computable metric space. Then a computable measure can be defined as a
computable point in $\mathcal{M}(X)$.

Let $\mathcal{X}=(X,d,S)$ be a computable metric space. Let us consider the
space $\mathcal{M}(X)$ of measures over $X$ endowed with weak topology, that
is: 
\begin{equation}
\mu_n\to\mu \mbox{ iff } \mu_n f \to \mu f 
\mbox{ for all real
continuous bounded } f
\end{equation}
where $\mu f$ stands for $\integral{f}{\mu}$.\newline

If $X$ is separable and complete, then $\mathcal{M}(X)$ is separable and
complete. Let $D\subset \mathcal{M}(X)$ be the set of those probability
measures that are concentrated in finitely many points of $S$ and assign
rational values to them. It can be shown that this is a dense subset (\cite%
{Bil68}). Let $\mu_{n_1,..,n_k,m_1,..,m_k}$ denote the measure concentrated
over the finite set $\{s_{n_1},\ldots,s_{n_k}\}$ with $q_{m_i}$ the weight
of $s_{n_i}$.

We consider Prokhorov metric $\rho$ on $\mathcal{M}(X)$ defined by: 
\begin{equation*}  \label{prokhorov}
\rho(\mu,\nu):=\inf \{\epsilon \in \mathbb{R}^+ :
\mu(A)\leq\nu(A^{\epsilon})+\epsilon \mbox{ for every Borel set }A\}.
\end{equation*}
where $A^{\epsilon}=\{x:d(x,A)< \epsilon\}$.\newline

This metric induces the weak topology on $\mathcal{M}(X)$. Furthermore, it
can be shown that the triple $(\mathcal{M}(X),D,\rho)$ is a computable
metric space (see \cite{Gac05}, \cite{HoyRoj07}).

\begin{definition}
\label{compmeas}A measure $\mu $ is computable if there is an algorithmic
enumeration of a fast sequence of ideal measures $(\mu _{n})_{n\in \mathbb{N}%
}\subset D$ converging to $\mu $ in the Prokhorov metric and hence, in the
weak topology.
\end{definition}

We need a criteria to check that a measure is computable. Let us then
introduce (following \cite{Gac05}) a certain fixed, enumerated sequence of
Lipschitz functions. Let $\mathcal{F}_0$ be the set of functions of the
form: 
\begin{equation}  \label{lip functions}
g_{s,r,\epsilon}=|1-|d(x,s)-r|^+/\epsilon|^+
\end{equation}
where $s\in S$, $r,\epsilon \in\mathbb{Q}$ and $|a|^+=\max\{a,0\}$.

These are Lipschitz functions equal to 1 in the ball $B(s,r)$, to 0 outside $%
B(s,r+\epsilon)$ and with intermediate values in between. It is easy to see
that the real valued functions $g_{s_i,r_j,\epsilon_k}:X\to \mathbb{R}$ are
computable, uniformly in $i,j,k$.


Let $\mathcal{F}$ be the smallest set of functions containing $\mathcal{F}_0$
and the constant 1, and closed under $\max$, $\min$ and rational linear
combinations. Clearly, this is also a uniform family of computable
functions. We fix some enumeration $\nu_{\mathcal{F}}$ of $\mathcal{F}$ and
we write $g_n$ for $\nu_{\mathcal{F}}(n)\in\mathcal{F}$. We remark that this
set is dense in the set of continuous functions with compact support.

The following lemma, proved in \cite{Gac05}, shows that the approach to
define computable measures we adopted, approximating measures with measures
supported on finite ideal sets is compatible with viewing the space of
measures as the dual of continuous functions, i.e. a measure is computable
if and only if it is a computable function $:C_{b}^{0}\rightarrow \mathbb{R}$%
.

\begin{lemma}
\label{gacs} Let $\mathcal{F}=\{g_1,g_2,...\}$ be the set introduced above.
A probability measure $\mu$ is computable if and only if $\integral{g_i}{\mu}$ is
computable uniformly in $i$.
\end{lemma}

Together with the previous lemma, the following result (see \cite{HoyRoj07})
will be all we use about computable measures:

\begin{lemma}
\label{lower_semi_compute} A probability measure $\mu$ is computable if and
only if the measure of finite union of ideal balls $\mu(B_{i_1}\cup\ldots%
\cup B_{i_k})$ is lower semi-computable, uniformly in $i_1,\ldots,i_k$.
\end{lemma}

\subsection{Computable probability spaces}
To obtain computability results on dynamical systems, it seems obvious that
some computability conditions must be required on the system. But the
``good'' conditions, if any, are not obvious to specify.

A computable function defined on the whole space is necessarily continuous.
But a transformation or an observable need not be continuous at every point,
as many interesting examples prove (piecewise-defined transformations,
characteristic functions of measurable sets,...), so the requirement of
being computable everywhere is too strong. In a measure-theoretical setting,
the natural weaker condition is to require the function to be computable 
\emph{almost everywhere}. It is not sufficient, and a computable condition
on the set on which the function is computable is needed to assure the
existence of computable points inside it. A condition which makes things
work is the following one:

\begin{definition}
A \defin{computable probability space} is a pair $(X,\mu)$ where $X$ is a
computable metric space and $\mu$ a computable Borel probability measure on
$X$.

Let $Y$ be a computable metric space. A function $(X,\mu)\to Y$ is \defin{$\delta$ a.e.-computable} if it is computable
on a constructive $G_\delta$-set of measure one, denoted by
$\dom f$ and called the \emph{domain of computability of $f$}.

A \defin{morphism} of computable probability spaces $f:(X,\mu)\to(Y,\nu)$ is a
morphism of probability spaces which is $\delta$ a.e.-computable.
\end{definition}

\begin{remark}\label{com_sum}
A sequence of functions $f_n$ is uniformly $\delta$ a.e.-computable if the functions are uniformly computable on their respective domains, which are uniformly constructive $G_\delta$-sets. Remark that in this case intersecting all the domains provides a constructive $G_\delta$-set on which all $f_n$ are computable. In the following we will apply this principle to the iterates $f_n=T^n$ of a  $\delta$ a.e.-computable function $T:X\to X$, which are uniformly  $\delta$ a.e.-computable.
\end{remark}


\section{Constructive Borel-Cantelli sets}

Given a space $X$ endowed with a probability measure $\mu ,$ the well known
Borel Cantelli lemma states that if a sequence of sets $A_{k}$ is such that
$\sum \mu (A_{k})<\infty $ then the set of points which belong to
finitely many $A_{k}$'s has full measure. In this section we show that if
the $A_{k}$ are given in some ``constructive'' way (and $\mu $ is computable) then this full
measure set contains some computable points. Hence this set contains points
which can effectively be constructed.

\begin{definition}
A sequence of positive numbers $a_{i}$ is \emph{effectively summable} if
the sequence of partial sums converges effectively: there is an algorithm
$A:\mathbb{Q\rightarrow N}$ such that if $A(\epsilon )=n$ then $\sum_{i\geq
n}a_{i}\leq \epsilon $.
\end{definition}

\begin{remark}
A computable sequence of \emph{positive }real numbers is effectively
summable if and only if its sum is a computable real number.
\end{remark}

For sake of simplicity, we will focus on the complements $U_n$ of the
$A_n$.

\begin{definition}
\label{ebc}A \emph{constructive Borel-Cantelli sequence} is a sequence 
$(U_{n})_{n\in \mathbb{N}}$ of uniformly constructive open sets such that the sequence
$(\mu(U_{n}^{\mathcal{C}}))$ is effectively summable.

The corresponding \emph{constructive Borel-Cantelli set} is $\bigcup_k
\bigcap_{n\geq k} U_n$.
\end{definition}

The Borel-Cantelli lemma says that every Borel-Cantelli set has
full-measure: we are going to see that every \emph{constructive} Borel-Cantelli
set contains a dense subset made of computable points.

\begin{lemma}[Normal form lemma]
\label{normal_form_lemma} Every constructive
Borel-Cantelli sequence can be effectively transformed into a constructive
Borel-Cantelli sequence $(U_{n})_{n\in \mathbb{N}}$ giving the same
Borel-Cantelli set, with $\mu (U_{n}^{\mathcal{C}})<2^{-n}$.
\end{lemma}

\begin{proof}
consider a constructive Borel-Cantelli sequence $(V_n)$. As
$\mu(X\setminus V_n)$ is effectively summable, an increasing sequence
$(n_i)_{i\geq 0}$ of integers can be computed such that for all $i,
\sum_{n\geq n_i} \mu(X\setminus V_n) < 2^{-i}$.

We now gather the $V_n$ by blocks, setting: 
\begin{equation*}
U_i := \bigcap_{n_i\leq n<n_{i+1}} V_n
\end{equation*}

$U_{i}$ is constructively open uniformly in $i$, and: 
\begin{equation*}
\mu (U_{i}^{\mathcal{C}})<2^{-i}\quad \mbox{ and }\quad
\bigcup_{k}\bigcap_{n\geq k}V_{n}=\bigcup_{i}\bigcap_{n\geq
n_{i}}V_{n}=\bigcup_{i}\bigcap_{j\geq i}U_{j}
\end{equation*}
\end{proof}

In the sequel we will always suppose that a constructive Borel-Cantelli
sequence is put in this normal form.

\begin{proposition}
Every finite intersection of constructive Borel-Cantelli sets is a constructive
Borel-Cantelli set.
\end{proposition}

\begin{proof}
let $(U_{n})$ and $(V_{n})$ be two constructive Borel-Cantelli sequences in
normal form. It is easy to see that:
\begin{equation*}
\bigcup_{k}\bigcap_{n\geq k}U_{n}\cap \bigcup_{k}\bigcap_{n\geq
k}V_{n}=\bigcup_{k}\bigcap_{n\geq k}U_{n}\cap V_{n}
\end{equation*}
and $\mu ((U_{n}\cap V_{n})^{\mathcal{C}})<2^{-n+1}$ which is effectively
summable.
\end{proof}

As every effectivity notion, the notion of constructive Borel-Cantelli set
naturally comes with its uniform version.

\begin{proposition}
\label{intersection_BC_proposition} The intersection of any uniform family
of constructive Borel-Cantelli sets contains a constructive Borel-Cantelli set.
\end{proposition}

\begin{proof}
suppose that $R_{i}=\bigcup_{k}\bigcap_{n\geq k}U_{n}^{i}$ is in normal
form. Consider a simple bijection $\varphi:\{(n,i):0\leq i\leq
n\}\rightarrow \mathbb{N}$ (for instance, $\varphi(n,i)=n(n+1)/2+i$)
computable in the two ways and define the sequences $(V_{m})_{m\in
\mathbb{N}}$ and $(a_{m})_{m\in \mathbb{N}}$ by $V_{m}=U_{n}^{i}$ and
$a_{m}=2^{-n}$ where $\varphi (n,i)=m$. Obviously
$\mu(V_{m}^{\mathcal{C}})>a_{m}$.

A simple calculation shows that $\sum a_m=4$. So $(V_m)$ is a constructive
Borel-Cantelli sequence.

Fix some $i$. If a point is outside $U_{n}^{i}$ for infinitely many $n$, it
is outside $V_{m}$ for infinitely many $m$. That is to say: 
\begin{equation*}
\bigcup_{k}\bigcap_{m\geq k}V_{m}\subseteq \bigcup_{k}\bigcap_{n\geq
k}U_{n}^{i}=R_{i}
\end{equation*}

As it is true for every $i$, the constructive Borel-Cantelli set induced by
$(V_m)_m$ is included in every $R_i$.
\end{proof}


\subsection{Computable points in constructive Borel-Cantelli sets}

The Borel-Cantelli lemma can be strengthened for constructive Borel-Cantelli
sequences obtaining that they contain computable points.

\begin{theorem}
\label{effective_BC_theorem} Let $X$ be a complete CMS and
$\mu$ a computable Borel probability measure on $X$.

For every constructive Borel-Cantelli set $R$, the set of computable points
lying in $R$ is dense in the support of $\mu$.
\end{theorem}

In order to the prove this theorem, we need the following lemma to
construct a computable point from what could be called a \emph{shrinking
sequence of constructive open sets}.

\begin{lemma}[Shrinking sequence]\label{lemma_shrinking}
Let $X$ be a complete CMS. Let $V_i$ be a sequence of non-empty
uniformly constructive open sets such that $\overline{V}_{i+1}\subseteq
V_i$ and $\diam(V_i)$ converges effectively to $0$. Then $\bigcap_i V_i$ is a
singleton containing a computable point.
\end{lemma}

\begin{proof}
As $V_i$ is non-empty there is a computable sequence of ideal points
$s_i\in V_i$. This is a Cauchy sequence, which converges by
completeness. Let $x$ be its limit: it is a computable point as
$\diam(V_i)$ converges to $0$ in an effective way. Fix some $i$: for all $j\geq
i$, $s_j\in V_j\subseteq \overline{V}_i$ so $x=\lim_{j\to\infty}s_j \in
\overline{V}_i$. Hence $x\in \bigcap_i \overline{V}_i= \bigcap_i V_i$.
\end{proof}

\begin{proof}[Proof of theorem \ref{effective_BC_theorem}]
Let $(U_n)_n$ be a constructive Borel-Cantelli sequence, in normal form ($\mu(U_n)>1-2^{-n}$, see lemma \ref{normal_form_lemma}). Let $B$ be an ideal
ball of radius $r\leq 1$ and positive measure. In $B$ we construct a computable
point which lies in $\bigcup_n\bigcap_{k\geq n}U_k$.

To do this, let $V_0=B$ and $n_0$ be such that $\mu(B)>2^{-n_0+1}$ (such an
$n_0$ can be effectively found from $B$): from this we construct a
sequence $(V_i)_i$ of uniformly constructive open sets and a computable increasing
sequence $(n_i)_i$ of natural numbers satisfying:
\begin{enumerate}
\item $\mu(V_i)+\mu(\bigcap_{k\geq n_i} U_k)>1$,
\item $V_i\subseteq \bigcap_{n_0\leq k<n_i} U_k$,
\item $\diam(V_i)\leq 2^{-i+1}$,
\item $\overline{V}_{i+1}\subseteq V_i$.
\end{enumerate}

The last two conditions assure that $\bigcap_i V_i$ is a computable point
(lemma \ref{lemma_shrinking}), the second condition assures that this point
lies in $\bigcap_{k\geq n_0} U_k$.

Suppose $V_i$ and $n_i$ have been constructed.

\begin{claim}
There exist $m>n_i$ and an ideal ball $B'$ of radius
$2^{-i-1}$ such that
\begin{equation}\label{claim}
\mu(V_i\cap\bigcap_{n_i\leq k<m}U_k\cap B')>2^{-m+1}.
\end{equation}
\end{claim}

We now prove of the claim:
By the first condition, $\mu(V_i\cap\bigcap_{k\geq n_i}U_k)>0$ so there
exists an ideal ball $B'$ of radius $2^{-i-1}$ such that
$\mu(V_i\cap\bigcap_{k\geq n_i}U_k\cap B')>0$. There is $m>n_i$ such that
$\mu(V_i\cap\bigcap_{k\geq n_i}U_k\cap B')>2^{-m+1}$, which implies the
assertion, and the claim is proved.

As inequality (\ref{claim}) can be semi-decided, such an $m$ and a $B'$ can
be effectively found. For $V_{i+1}$, take any finite union of balls whose
closure is contained in $V_i\cap\bigcap_{n_i\leq k<m}U_k\cap B'$ and whose measure is greater than
$2^{-m+1}$. Put $n_{i+1}=m$. Conditions 2., 3. and 4. directly follow from
the construction, condition 1. follows from
$\mu(V_{i+1})>2^{-m+1}>1-\mu(\bigcap_{k\geq m}U_k)$ (the sequence is in normal form).
\end{proof}

The following corollary allows to apply the above criteria to a uniform
infinite sequence of constructive Borel-Cantelli sets.

\begin{corollary}
\label{intersection_BC_corollary} Let $X$ be a complete CMS  and $\mu$ a
computable Borel probability measure on $X$.

Given a uniform family $(R_{i})_{i}$ of constructive Borel-Cantelli sets, the
set of computable points lying in $\bigcap_{i}R_{i}$ is dense in the support
of $\mu$.
\end{corollary}

\begin{proof}
this a direct consequence of proposition \ref{intersection_BC_proposition}
and theorem \ref{effective_BC_theorem}.
\end{proof}

We remark that, in the particular case of Cantor spaces with an uniform
measure a result of this kind can also be obtained from \cite{Sch71} since it
is possible to relate Borel Cantelli sequences to Schnorr tests. This
relation is developped in \cite{GacHoyRoj08} giving some new connections between
Schnorr randomness and dynamical typicality.

\subsubsection{Application to convergence of random variables}
Here, $(X,\mu)$ is a computable probability space.

\begin{definition}
A \defin{random variable} on $(X,\mu)$ is a measurable function
$f:X\to \R$.
\end{definition}

\begin{definition}
Random variables $f_n$ \defin{effectively converge in probability} to $f$ if for each $\epsilon>0$, $\mu\{x:
|f_n(x)-f(x)|<\epsilon\}$ converges effectively to $1$, uniformly in
$\epsilon$. That is, there is a computable function $n(\epsilon,\delta)$
such that for all $n\geq n(\epsilon,\delta)$, $\mu[|f_n-f|\geq\epsilon]<\delta$.
\end{definition}

\begin{definition}
Random variables $f_n$ \defin{effectively converge almost surely} to $f$ if
$\sup_{k\geq n}|f_n-f|$ effectively converge in probability to $0$.
\end{definition}




\begin{theorem}\label{theorem_convergence_Borel_Cantelli}
Let $f_n,f$ be uniformly  $\delta$ a.e.-computable random variables. If $f_n$ effectively
converges almost surely to $f$ then the set $\{x:f_n(x)\to f(x)\}$ contains
a constructive Borel-Cantelli set.

In particular, the set of computable
points for which the convergence holds is dense in $\supp(\mu)$.
\end{theorem}

\begin{proof}
Let $D=\bigcap_n D_n$ be a constructive $G_\delta$-set of full measure on which all $f_n,f$ are computable. $D_n$ are uniformly constructive open sets, and we can suppose $D_{n+1}\subseteq D_n$ (otherwise, replace $D_n$ by $D_0\cap\ldots\cap D_n$).

There are uniformly constructive open sets
$U_n(\epsilon)$ such that $U_n(\epsilon)\cap D=[|f_n-f|<\epsilon]\cap D$. $\mu(\bigcap_{n\geq k}U_n(\epsilon))$
converges effectively to $1$, uniformly in $\epsilon$ so it is possible to
compute a sequence $(k_i)_i$ such that $\mu(\bigcap_{n\geq
k_i}U_n(2^{-i}))>1-2^{-i}$ for all $i$. Put $V_i=\bigcap_{k_i\leq
n<k_{i+1}}U_n(2^{-i})\cap D_i$: $V_i$ is constructively open uniformly in $i$ and
$\mu(V_i)>1-2^{-i}$. The sets $V_i$ form a constructive Borel-Cantelli sequence,
and if a point $x$ is in the corresponding Borel-Cantelli set then $x\in D$ and there is
$i_0$ such that $x\in V_i$ for all $i\geq i_0$, so $|f_n(x)-f(x)|<2^{-i}$
for all $n\geq k_i,i\geq i_0$. Hence $f_n(x)\to f(x)$.
\end{proof}


\section{Pseudorandom points and dynamical systems}

Let $X$ be a metric space, let $T:X\mapsto X$ be a Borel map. Let $\mu $ be
an invariant Borel measure on $X$, that is: $\mu (A)=\mu (T^{-1}(A))$ holds
for each measurable set $A$. A set a $A$ is called $T$-invariant if
$T^{-1}(A)=A(\text{mod } 0)$. The system $(T,\mu )$ is said to be ergodic
if each $T$-invariant set has total or null measure. In such systems the
famous Birkhoff ergodic theorem says that time averages computed along $\mu
$ typical orbits coincides with space average with respect to $\mu .$ More
precisely, for any $f\in L^{1}(X)$ and it holds
\begin{equation}
\underset{n\rightarrow \infty }{\lim }\frac{S_{n}^{f}(x)}{n}=\integral{f}{\mu},
\label{Birkhoff}
\end{equation}
for $\mu $ almost each $x$, where $S_{n}^{f}=f+f\circ T+\ldots +f\circ
T^{n-1}.$

If a point $x$ satisfies equation \ref{Birkhoff} for a certain $f$, then we
say that $x$ is typical with respect to the observable $f$.

\begin{definition}
\label{mutyp}If $x$ is typical w.r.t any continuous function $f:X\to \mathbb{R}$ with compact support, then we call it a $\boldsymbol{\mu}$\textbf{\emph{-typical
point}}.
\end{definition}

In this section we will see how the constructive Borel-Cantelli lemma can be
used to prove that in a large class of interesting systems there exists
computable typical points.

Let us call $(X,\mu,T)$ a \textbf{\emph{computable ergodic system}} if $(X,\mu)$ is a
computable probability space, $T$ is an endomorphism (i.e. an  $\delta$ a.e.-computable measure-preserving
transformation) and $(X,\mu,T)$ is ergodic.

Before to enter in the main theme of typical statistical behaviors let us
see an easier topological result in this line. One of the features of
undecomposable (topologically transitive) chaotic systems is that there are
many dense orbits, the following shows that if the system is computable then
there are computable dense orbits.

We remark that this result can also be obtained as a corollary of the
constructive Baire theorem \cite{YasMorTsu99}.

\begin{theorem}
\label{computable_dense} Let $\mathcal{X}$ be a computable complete metric
space and $T:X\rightarrow X$ a transformation which is computable on a dense
constructive open set. If $T$ has a dense orbit, then it has a computable one which
is dense.
\end{theorem}

In other words, there is a computable point $x\in X$ whose orbit is dense in 
$X$. Actually, the proof is an algorithm which takes an ideal ball as input
and computes a transitive point lying in this ball.

\begin{proof}
$(B_{i})_{i\in \mathbb{N}}$ being an enumeration of all ideal balls, define the open sets
$U_{i}=dom(f)\cap \bigcup_{n}T^{-n}B_{i}$ which are constructive uniformly in $i$. By
hypothesis, $U_{i}$ is also dense. $\bigcap_{i}U_{i}$ is the set of
transitive points. From any ideal ball $B(s_{0},r_{0})$ we effectively
construct a computable point in $B(s_{0},r_{0})\cap \bigcap_{i}U_{i}$.

If $B(s_i,r_i)$ has been constructed, as $U_i$ is dense $B(s_i,r_i)\cap
U_i$ is a non-empty constructive open set, so an ideal ball $B(s,r)\subseteq
B(s_i,r_i)\cap U_i$ can be effectively found (any of them can be chosen,
for instance the first coming in the enumeration). We then
set $B(s_{i+1},r_{i+1}):=B(s,r/2)$.

The sequence of balls computed satisfies: 
\begin{equation*}
\overline{B}(s_{i+1},r_{i+1}) \subseteq B(s_i,r_i) \cap U_0\cap\ldots\cap U_i
\end{equation*}

As $(r_{i})_{i\in \mathbb{N}}$ is a decreasing computable sequence
converging to $0$ and the space is complete, $(s_{i})_{i\in \mathbb{N}}$
converges effectively to a computable point $x$. Then
$\{x\}=\bigcap_{i}B(s_{i},r_{i})\subseteq \bigcap_{i}U_{i}$.
\end{proof}


\subsection{Computable typical points\label{comptypsec}}

We will use the results from the previous section to prove that computable
typical points exist for a class of dynamical systems. Each time the set of
typical points is a constructive Borel-Cantelli set, theorem \ref{typical_f}
applies.

For instance, in the case of the shift on the Cantor space with a Bernoulli
measure, the Birkhoff ergodic theorem reduces to the strong law of large
numbers, which proof is simpler and makes explicit use of the Borel-Canteli
lemma. This is possible thanks to the independence between the random
variables involved, but strict independence is actually unnecessary: the
proof can be adapted whenever the correlations between the random variables
decrease sufficiently fast.

\begin{definition}
\label{logarithmically ergodic}We say that a system $(X,T,\mu )$ is
\textbf{\emph{$\boldsymbol{\ln^2}$-ergodic}} for observables in some set of
functions $\mathcal{B}$ if for each $(\phi,\psi)\in \mathcal{B}^{2}$ there
is $c_{\phi,\psi}>0$ such that
\begin{equation*}
\left|\frac{1}{n}\sum_{i<n}\integral{\phi\circ
T^i\psi}{\mu}-\integral{\phi}{\mu}\integral{\psi}{\mu}\right|\leq
\frac{c_{\phi ,\psi }}{(\ln(n))^2}\qquad \text{ for all }n\geq 2.
\end{equation*}
\end{definition}

Now we can state:
\begin{theorem}\label{theorem_effective_convergence} Let $(X,T,\mu)$ be a dynamical system which is $\ln^2$-ergodic for observables in some set $\mathcal{B}$ of bounded observables. For each $\phi\in\mathcal{B}$, the almost-sure convergence:
\[
\frac{1}{n}\sum_{i<n}\phi\circ T^i\to_n\integral{\phi}{\mu}
\]
is effective.
\end{theorem}

Note that for the moment, no computability assumption is needed on the system.

As announced, the proof is an adaptation of the proof of the strong law of
large numbers. We first prove two lemmas.

\begin{lemma}\label{lemma_sequence}
There exists a computable sequence $n_i$ such that:
\begin{itemize}
\item $\beta_i:=\frac{n_i}{n_{i+1}}$ converge effectively to $1$,
\item $\frac{1}{\ln(n_i)^2}$ is effectively summable.
\end{itemize}
\end{lemma}

\begin{proof}
For instance, take $n_i=\ceil{(1+i^{-\alpha})^i}$ with $0<\alpha<1/2$.
\end{proof}

From now on, we denote $\frac{S^\phi_n}{n}$ by $f_n$. 

\begin{lemma}\label{lemma_subsequence}
The almost-sure convergence of the subsequence $f_{n_i}$ to $\integral{\phi}{\mu}$ is effective.
\end{lemma}

\begin{proof}
For $\delta >0$, define the deviation sets: 
\begin{equation*}
A_{n}^{\phi }(\delta )=\left\{ x\in X:\left\vert f_n(x)-\integral{\phi}{\mu} \right\vert \geq \delta \right\} .
\end{equation*}

By Tchebytchev inequality, 
\begin{equation*}
\delta^2\mu (A_{n}^{\phi}(\delta ))\leq \norm{f_n-\integral{\phi}{\mu}}{L^{2}}^2.
\end{equation*}

Since adding a constant to $\phi $ does not change this quantity, without
loss of generality, let us suppose that $\integral{\phi}{\mu}=0$. Then 
\begin{equation*}
\norm{\frac{S_{n}^{\phi}}{n}-\integral{\phi}{\mu}}{L^{2}}^{2}=\integral{\left(\frac{S_{n}^{\phi}}{n}\right)^{2}}{\mu}
=\integral{\left(\frac{\phi +\phi \circ T+...+\phi \circ
T^{n-1}}{n}\right)^{2}}{\mu}
\end{equation*}

by invariance of $\mu $ this is equal to 
\begin{equation*}
\frac{1}{n^{2}}\integral{n\phi^{2}}{\mu}
+\frac{2}{n^{2}}\integral{\big(\underset{i<j<n}\sum\phi\circ
T^{j-i}\phi\big)}{\mu}
\end{equation*}

hence,

\begin{eqnarray*}
\delta^2 \mu (A_{n}^{\phi}(\delta )) & \leq & \frac{M^2}{n}+\frac{2}{n}\sum_{k<n}\integral{\phi\circ T^k\phi}{\mu} \\
& \leq & \frac{M^2}{n}+\frac{c_{\phi,\phi}}{\ln(n)^2}
\end{eqnarray*}

As $\frac{M^2}{n_i}+\frac{c_{\phi,\phi}}{\ln(n_i)^2}$ is effectively summable (by choice of $n_i$, see lemma \ref{lemma_sequence}) uniformly in $\delta$, it follows that $f_{n_i}$ converge effectively almost-surely to $\integral{\phi}{\mu}$.
\end{proof}

As $n_i$ is not dispersed too much, the almost-sure convergence of the subsequence $f_{n_i}$ implies that of the whole sequence $f_n$. Actually the effectivity is also preserved. We now make this precise.

\begin{lemma}
For $n_i\leq n<n_{i+1}$ and $\beta_{i}:=\frac{n_i}{n_{i+1}}$, one has:

\begin{equation}\label{variation}
\left\Vert f_{n_i}-f_n\right\Vert_\infty \leq 2(1-\beta_{i})\left\Vert \phi \right\Vert_\infty.
\end{equation}
\end{lemma}

\begin{proof}
Let $M=\norm{\phi}{\infty}$. To see this, for any $k,l,\beta$ with $\beta \leq k/l \leq 1$: 
\begin{eqnarray*}
\frac{S_{k}^{\phi}}{k}-\frac{S_{l}^{\phi}}{l} &=&\left( 1-\frac{k}{l}\right) 
\frac{S_{k}^{\phi}}{k}-\frac{S_{l-k}^{\phi}\circ T^{l-k}}{l} \\
&\leq &(1-\beta )M+\frac{(l-k)M}{l}=2(1-\beta)M,
\end{eqnarray*}
Taking $\beta =\beta_{i}$ and $k=n_{i}$, $l=n$ first and then $k=n$,
$l=n_{i+1}$ gives the result.
\end{proof}

\begin{proof}[Proof of theorem \ref{theorem_effective_convergence}]
Let $\delta,\epsilon>0$. To prove that $f_n$ converge effectively almost-surely, one has to compute some $p$ (from $\delta$ and $\epsilon$) such that $\mu(\bigcup_{n\geq p}A_n(\delta))<\epsilon$.

As $\beta_i$ converge effectively to $1$, one can compute $i_0$ such that if $i\geq i_0$ then $\beta_i>1-\delta/(4M)$. Inequality \ref{variation} then implies
\[
\bigcup_{n_i\leq n<n_{i+1}}A_n(\delta)\subseteq A_{n_i}(\delta/2).
\]
Indeed if $n_i\leq n<n_{i+1}$ and $|f_{n_i}(x)-\integral{\phi}{\mu}|<\delta/2$ then $|f_n(x)-\integral{\phi}{\mu}|\leq |f_n(x)-f_{n_i}(x)|+|f_{n_i}(x)-\integral{\phi}{\mu}|\leq \delta$.

As $f_{n_i}$ converge effectively almost-surely, one can compute some $j_0$ such that $\mu(\bigcup_{j\geq j_0} A_{n_j}(\delta/2))<\epsilon$. Let $p=n_k$ where $k=\max(i_0,j_0)$: $\bigcup_{n\geq p}A_n(\delta)\subseteq \bigcup_{j\geq j_0}A_{n_i}(\delta/2)$ whose measure is less than $\epsilon$.
\end{proof}

\begin{corollary}
\label{typical_f} Let $(X,T,\mu)$ be a \emph{computable} dynamical system which is $\ln^2$-ergodic for observables in some set
$\mathcal{B}$ of bounded functions and let $\phi$ be a  $\delta$ a.e.-computable
observable in $\mathcal{B}$.

The set of points which are typical w.r.t
$\phi$ contains a constructive Borel-Cantelli set. In particular, it contains computable  points.
\end{corollary}

\begin{proof}
Apply theorem \ref{theorem_convergence_Borel_Cantelli} to the sequence of uniformly  $\delta$ a.e.-computable functions $f_n=\frac{S^\phi_n}{n}$ which converge effectively almost-surely by theorem \ref{theorem_effective_convergence}.
\end{proof}

\begin{remark}
\label{renun}In the proof of thm \ref{theorem_effective_convergence} we see that the constructive Borel-Cantelli set
depends in an effective way on $\norm{\phi}{\infty}$ and $c_{\phi ,\phi }$. This gives
the possibility to operate in a way to apply
Prop. \ref{intersection_BC_proposition} and
Cor. \ref{intersection_BC_corollary} to find a constructive Borel Cantelli
set and computable points contained in the set of points typical with
respect to a uniform family $\phi_{i},T_{i}$.
\end{remark}

By the above remark, to construct $\mu $-typical points (see definition
\ref{mutyp}) using the above mentioned results, the following conditions
are sufficient:

\begin{theorem}
\label{mutypicalteo}If a computable system is $\ln^2$-ergodic for observables in $\mathcal{F}=\{g_{1},g_{2},\ldots \}$ (this
set was defined in section \ref{seccompmu}) and the associated constants
$c_{g_{i}}$ (see definition \ref{polynomial decay}) can be estimated
uniformly in $i$ (there is an algorithm $A:\mathbb{N\rightarrow Q}$ such
that $A(i)\geq c_{g_{i}}$) then it has a set of computable $\mu$-typical
points which is dense in the support of $\mu$.
\end{theorem}

\begin{proof}
We remark that $\mathcal{F}$ is dense in the set of continuous functions on
$X$ with compact support (with the sup norm) hence a computable point
which is typical for each $g_{i}$ is $\mu$-typical. Such points can be
found by applying theorem \ref{typical_f} for each $g_{i}$ and using
proposition \ref{intersection_BC_proposition} as explained in remark
\ref{renun}.
\end{proof}

\subsubsection{$\ln^2$-mixing}
We will apply this to systems having a stronger property: they are
\emph{mixing}, with logarithmical speed. More precisely, this can be quantified using
the \textbf{\emph{correlation functions}}:

\begin{equation*}
C_n(\phi,\psi)=\left|\integral{\phi \circ T^{n}\psi}{\mu} -\integral{\phi}{\mu} \integral{\psi}{\mu}\right|
\end{equation*}

which measures the dependence between observation through $\phi$ and $\psi$
at times $n \gg 1$ and $0$ respectively (possibly with $\psi=\phi$). Note
that $C_n(\phi,\psi)=0$ corresponds, in probabilistic terms, to $\phi\circ
T^n$ and $\psi$ being independent random variables.

\begin{definition}
\label{polynomial decay}We say that a system $(X,T,\mu )$ has
\textbf{\emph{$\boldsymbol{\ln^2}$-decay of correlations}} for observables in some
set of functions $\mathcal{B}$ if for each $(\phi
,\psi )\in \mathcal{B}^{2}$ there is $c_{\phi ,\psi }>0$ such that
\begin{equation*}
C_{n}(\phi ,\psi )\leq \frac{c_{\phi ,\psi }}{(\ln(n))^2}\qquad \text{ for
all }n\geq 2.
\end{equation*}
\end{definition}

\begin{lemma}
If a system has $\ln^2$-decay of correlation for observables in
$\mathcal{B}$ then it is $\ln^2$-ergodic for observables in $\mathcal{B}$. The ergodicity constants depend in an effective way on the
mixing constants.
\end{lemma}

\begin{proof}
We first prove that for all $n\geq 2$,
\begin{equation}\label{eq1}
\sum_{k=2}^n \frac{1}{\ln(k)^2} \leq \frac{2n}{\ln(n)^2}+4
\end{equation}

For $n\geq 56$,
\begin{eqnarray*}
\sum_{k=56}^n \frac{1}{\ln(k)^2} & \leq & \int_{x=55}^n\!\frac{\mathrm{d}x}{\ln(x)^2} \\
& \leq &
\int_{x=55}^n\!2\left(\frac{1}{\ln(x)^2}-\frac{2}{\ln(x)^3}\right)\mathrm{d}x
\quad \text{ (as $55\geq \ln(4)$)}\\
& = & \frac{2n}{\ln(n)^2}-\frac{110}{\ln(55)^2}
\end{eqnarray*}
which, combined with $\sum_{k=2}^{55} \frac{1}{\ln(k)^2}\leq 10$ and
$\frac{110}{\ln(55)^2}\geq 6$, gives inequality (\ref{eq1}).

Finally, for $n\geq 2$,
\begin{eqnarray*}
\left|\frac{1}{n}\sum_{i<n}\integral{\phi\circ
T^i\psi}{\mu}-\integral{\phi}{\mu}\integral{\psi}{\mu}\right| & \leq &
\frac{1}{n}\sum_{i<n}C_i(\phi,\psi) \\
& \leq & \frac{2c_{\phi,\psi}}{\ln(n)^2} + \frac{4c_{\phi,\psi}}{n}\\
& \leq & \frac{6 c_{\phi,\psi}}{\ln(n)^2}
\end{eqnarray*}
\end{proof}

\subsection{Application: computable absolutely normal numbers\label{normal}}

An absolutely normal (or just normal) number is, roughly speaking, a real
number whose digits (in every base) show a uniform distribution, with all
digits being equally likely, all pairs of digits equally likely, all
triplets of digits equally likely, etc.

While a general, probabilistic proof can be given that \emph{almost all}
numbers are normal, this proof is not constructive and only very few
concrete numbers have been shown to be normal. It is for instance widely
believed that the numbers $\sqrt{2}$, $\pi $ and $e$ are normal, but a proof
remains elusive. The first example of an absolutely normal number was given
by Sierpinski in 1916, twenty years before the concept of computability was
formalized. Its construction is quite complicate and is a priori unclear
whether his number is computable or not. In \cite{BechFig02} a recursive
reformulation of Sierpinski's construction (equally complicate) was given,
furnishing a computable absolutely normal number.

As an application of theorem \ref{typical_f} we give a simple proof that
computable absolutely normal numbers are dense in $[0,1]$.

Let $b$ be an integer $\geq 2$, and $X_b$ the space of infinite sequences
on the alphabet $\Sigma_{b}=\{0,\ldots ,b-1\}$. Let $T=\sigma $ be the
shift transformation on $X_b$, and $\lambda$ be the uniform measure. A real
number $r\in[0,1]$ is said to be absolutely normal if for all $b\geq 2$,
its $b$-ary expansion $r_b\in X_b$ satisfies:
\begin{equation*}
\lim_{n\rightarrow \infty} \frac{1}{n}\sum_{i=0}^{n-1}1_{[w]}\circ
\sigma^i(r_b)=\frac{1}{b^{|w|}} \qquad \mbox{ for all } w\in \Sigma_b^*.
\end{equation*}

\begin{theorem}
\label{compu_normal}The set of computable reals which are absolutely normal
is dense in $[0,1]$.
\end{theorem}

\begin{proof}
for each base $b\geq 2$, consider the transformation
$T_{b}:[0,1]\rightarrow \lbrack 0,1]$ defined by $T_{b}(x)=bx(\mbox{mod
}1)$. The Lebesgue measure $\lambda $ is $T_{b}$-invariant and
ergodic. The partition in intervals $[k/b,(k+1)/b[$ induces the symbolic
model $(\Sigma _{b}^{\mathcal{N}},\sigma ,\lambda )$ which is
measure-theoretically isomorphic to $([0,1],T_{b},\lambda )$: the interval
$[k/b,(k+1)/b[$ is represented by $k\in \Sigma _{b}$. For any word $w\in
\Sigma _{b}$ define $I(w)$ to be the corresponding interval
$[0.\overline{w},0.\overline{w}+2^{-|w|}]$.

Defining $\dom T_{b}:=[0,1]\setminus \{\frac{k}{b}:0\leq k\leq b\}$
(the interior of the partition) makes $T_{b}$ a  $\delta$ a.e.-computable
transformation. The observable $f_{w}:=1_{I(w)}$ is also
 $\delta$ a.e.-computable, with $\dom f_{w}=[0,1]\setminus \partial I(w)$.

Actually, since $f_w\circ \sigma^n$ and $f_w$ are independent for $n>|w|$,
theorem \ref{typical_f} applies to $([0,1],T_{b},\lambda )$ and $f_{w}$.
Therefore, the set of points (for the system $(T_b,\lambda)$) which are
typical w.r.t the observable $f_w$ contains a constructive Borel-Cantelli set
$R_{b,w}$. Furthermore, $R_{b,w}$ is constructive \emph{uniformly} in $b,w\in
\Sigma_b$. Hence, by corollary \ref{intersection_BC_corollary}, their
intersection, which is made of absolutely normal numbers, contains a dense
set of computable points.
\end{proof}

\section{Dynamical systems having computable typical points\label{examples}}

We will see that in a large class of dynamical systems which have a single 
\emph{physically relevant} invariant measure, the computability of this
measure and related $c_{g_{i}}$, for observables in $\mathcal{F}$ can be
proved, hence we can apply Thm. \ref{mutypicalteo} to find pseudorandom
points in such systems.

\subsection{Physical measures}

In general, given $(X,T)$ there could be infinitely many invariant measures
(this is true even if we restrict to probability measures). Among this class
of measures, some of them are particularly important. Suppose that we
observe the behavior of the system $(X,T)$ trough a class of continuous
functions $f_{i}:X\rightarrow \mathbb{R}$. We are interested in the
statistical behavior of $f_{i}$ along typical orbits of the system. Let us
suppose that the time average along the orbit of $x$ exists
\begin{equation*}
A_{x}(f_{i})=\underset{n\rightarrow \infty }{\lim }\frac{1}{n}\sum
f_{i}(T^{n}(x))
\end{equation*}
this is a real number for each $f_{i}.$ Moreover $A_{x}(f_{i})$ is linear
and continuous with respect to small changes of $f_{i}$ in the $\sup $ norm.
Then the orbit of $x$ acts as a measure $\mu _{x}$ and $A_{x}(f_{i})=\integral{f_i}{\mu}_{x}$ (moreover this measure is also invariant for $T$). This
measure is physically interesting if it is given by a \textquotedblleft
large\textquotedblright\ set of initial conditions. This set will be called
the basin of the measure. If $X$ is a manifold, it is said that an invariant
measure is \emph{physical} (or SRB from the names of Sinai, Ruelle and
Bowen) if its basin has positive Lebesgue measure (see \cite{Y} for a survey
and more precise definitions).

In what follows we will consider SRB measures in the classes of systems
listed below,

\begin{enumerate}
\item The class of uniformly hyperbolic system on submanifolds of
$\mathbb{R}^n$.

\item The class of piecewise expanding maps on the interval.

\item The class of \emph{Manneville-Pomeau} type maps (non uniformly
expanding with an indifferent fixed point).
\end{enumerate}

All these systems, which are rather well understood, have a unique physical
measure with respect to which correlations decay is at least polynomial.
Furthermore, in each case, the corresponding constants can be estimated for
functions in $\mathcal{F}$. The computability of the physical measures is
proved case by case, but it is always a consequence of the fact that, in one
way or another, the physical measure is \textquotedblleft
approached\textquotedblright\ by iterates of the Lebesgue measure at a known
speed.

\subsection{Uniformly hyperbolic systems}

To talk about SRB measures on a system whose phase space is a manifold, we
have to introduce the Lebesgue measure on a manifold and check that it is
computable.

\subsubsection{Computable manifolds and the Lebesgue measure}

For simplicity we will not consider general manifolds but submanifolds of
$\mathbb{R}^{n}$.

\begin{definition}
Let $M$ be a computable metric subspace of $\mathbb{R}^{n}$. We say that
$M$ is a $m$ dimensional computable $C^{k}$ submanifold of $\mathbb{R}^{n}$
if there exists a computable function $f:M\times B(0,1)\rightarrow M$
(where $B(0,1)$ is the unit ball of $\mathbb{R}^{m}$and $M\times B(0,1)$
with the euclidean distance is a CMS in a natural way) such that for each
$x\in M$, $f_{x}=f(x,.)$ is a $C^{k}$ diffeomorphism with all $k$
derivatives being computable.
\end{definition}

For each $x$, the above $f_x$ is a map whose differential at any $z\in
B(\mathbf{0},1)$ is a linear, rank $m$ function
$Df_{x,z}:\mathbb{R}^{m}\rightarrow \mathbb{R}^{n}$. This can be seen as a
composition of two functions $Df_{x,z}=Df^2_{x,z}\circ Df^1_{x,z}$ such
that $Df^1_{x,z}:\mathbb{R}^{m}\rightarrow \mathbb{R}^{m}$ is invertible
and $Df^2_{x,z}:\mathbb{R}^{m}\rightarrow \mathbb{R}^{n}$ is an isometry.

Let us denote $B_{x}$ the image of $B(0,1)$ by $f_{x}$. Then the Lebesgue
measure of $D\subset B_{x}$ is defined as
\begin{equation*}
m(D)=\int_{f_{x}^{-1}(D)}\!\det (Df_{x,z}^{1})\,\mathrm{d}z.
\end{equation*}
This does not depend on the choice of $B_{x}$ and $f_{x}$, and it give rise
to a finite measure (Lebesgue measure) on $M$ (see \cite{GMS} page 74). This
measure is indeed the $m$ dimensional Hausdorff measure on $M.$ Moreover,
the Lebesgue measure is a computable measure.

\begin{lemma}
The Lebesgue measure on a computable $C^k$ submanifold of $\mathbb{R}^n$ is
computable.
\end{lemma}

\begin{proof}[Proof (sketch)]
Suppose that $A$ is a constructive open subset of some $B_s$, where $s$ is an ideal
point of $M$. Since the function $\det(Df^1_{x,z})$ is computable and the
function $1_{f^{-1}(A)}(z)$ is lower semi-computable, we can lower
semi-compute the value $m(A)$. In particular, there is a base of ideal
balls whose measures are lower semi-computable. Let $B$ and $B^{\prime}$ be
such balls. Since these balls have zero measure boundaries, we can compute
the measure of their intersection (which is a constructive open included in
$B$). Hence any constructive open set can be decomposed into a (same measure)
disjoint union of constructive open sets whose measures can be lower
semi-computed. By lemma \ref{lower_semi_compute}, $m$ is computable.
\end{proof}

\subsubsection{The SRB measure of uniformly hyperbolic systems}

Let us consider a connected $C^{2}$ computable manifold $M.$ Let us consider
a dynamical system $(M,T)$ where $T$ is a $C^{2}$ computable diffeomorphism
on $M.$

Let us consider a constructive open forward invariant set $Q\subset M$
(i.e. $T(\overline{Q})\subset Q$). Let us consider the (attracting) set
\begin{equation*}
\Lambda =\underset{n\geq 0}{\cap }T^{n}(Q).
\end{equation*}
Suppose that $\Lambda $ contains a dense orbit and that it is an hyperbolic
set for $T$, which means that the following conditions are satisfied.

There is a splitting of the tangent bundle of $M$ on $\Lambda $: $T_{\Lambda
}M=E_{\Lambda }^{s}\oplus E_{\Lambda }^{u}$ (at each point $x$ of $\Lambda $
the tangent space at $x$ can be splitted in a direct sum of two spaces, the
stable directions and the unstable ones) and a $\lambda _{0}<1$ such that

\begin{itemize}
\item the splitting is compatible with $T,$ that is:
$DT_{x}(E_{x}^{s})=E_{T(x)}^{s}$ and
$DT_{x}^{-1}(E_{x}^{u})=E_{T^{-1}(x)}^{u}$.

\item The dynamics expand exponentially fast in the unstable directions and
contracts exponentially fast in the stable directions in an uniform way,
that is: for each $x\in \Lambda $ and for each $v\in E_{x}^{s}$ and$\ w\in
E_{x}^{u},$ $\abs{DT_{x}(v)}\leq \abs{\lambda _{0}v}$ and $\abs{DT_{x}^{-1}(w)}\leq
\abs{\lambda _{0}w}$.
\end{itemize}

Under these assumptions it is known that

\begin{theorem}
(see \cite{Via97} e.g.)There is a unique invariant SRB measure $\mu $
supported on $\Lambda $. Moreover the measure is ergodic and its basin has
full Lebesgue measure on $Q$.
\end{theorem}

This measure has many good properties: it has exponential decay of
correlations and it is stable under perturbations of $T$ (see \cite{Via97}
e.g.). Another good property of this measure is that it is computable.

\begin{theorem}
\label{unifhyp}If $M$ and $T$ are $C^{2}$, computable and uniformly
hyperbolic as above, then the SRB measure $\mu $ is computable.
\end{theorem}

\begin{proof}
Let $m$ be the Lebesgue measure on $Q$ normalized by $m(Q)=1$, clearly it is
a computable measure. From \cite{Via97} (Prop. 4.9, Remark 4.2) it holds
that there are $\lambda <1$ such that for each $\nu $-H\"{o}lder ($\nu \in
(0,1]$) continuous observable $\psi $, it holds 
\begin{equation*}
\abs{\integral{\psi \circ T^{n}}{m}-\integral{\psi}{\mu}}\leq \lambda ^{n}~c_{\psi }
\end{equation*}
where $c_{\psi }=C\integral{\abs{\psi}}{m}+\norm{\psi}{\nu }$, where $C$
is independent from $\psi $ and then can be estimated for each uniform
sequence $\psi _{i}\in \mathcal{F}$ uniformly in $i$. This means that for
each $\psi _{i}\in \mathcal{F}$ its integral with respect to $\mu $ can be
calculated up to any given accuracy, uniformly in $i$. Indeed if we want to
calculate $\integral{\psi_i}{\mu}$ up to an error of $\epsilon $ we
calculate $c_{\psi _{i}}$ up to an error of $\epsilon $ (this error is not
really important as we will see immediately) and choose an $n$ such that
$c_{\psi _{i}}\lambda ^{n}\leq \frac{\epsilon }{2}.$

By this we know that $\abs{\integral{\psi_i\circ
T^{n}}{m}-\integral{\psi_i}{\mu}}\leq \frac{\epsilon }{2}$. Now we have to
calculate $\integral{\psi_i\circ T^{n}}{m}$ up to an error of
$\frac{\epsilon }{2}$ and this will be the output. By lemma \ref{gacs} then
$\mu $ is computable.
\end{proof}

\begin{corollary}
In an unif. hyp. computable system equipped with its SRB measure as above,
the set of computable $\mu -$typical points is dense in the support of $\mu$.
\end{corollary}

\begin{proof}
$\mu $ is computable by the previous theorem, and the correlations decay is
given by proposition 4.9 in \cite{Via97} from which follows that there is
$\lambda <1$ such that for each $(g_{i},g_{j})\in \mathcal{F}^{2}$ it
holds,
\begin{equation*}
\abs{\integral{g_{i}\circ T^{n}g_{j}}{m}-\integral{g_i}{\mu}
\integral{g_j}{\mu}}\leq \lambda ^{n}~c_{g_{i},g_{j}}
\end{equation*}
where $c_{g_{i},g_{j}}=C(\integral{|g_i|}{m}+\norm{g_{i}}{1})(\integral{|g_{j}|}{m}+\norm{g_{j}}{1})$ ($C$ is a constant
independent of $g_{i}\in \mathcal{F}$, $\norm{\ast}{1}$ is the
Lipschitz norm, since functions in $\mathcal{F}$ are Lipschitz) are
computable uniformly in $i,j$.  Then the result follows from theorem
\ref{mutypicalteo}.
\end{proof}

\subsection{Piecewise expanding maps}

We introduce a class of discontinuous maps on the interval having an
absolutely continuous SRB invariant measure. The density of this measure has
also bounded variation. We will show that this invariant measure is
computable.

Let $I$ be the unit interval. Let $T:I\rightarrow I$ we say that $T$ is
piecewise expanding if there is a finite partition $P=\{I_{1},...,I_{k}\}$
of $I,$ such that $I_{i}$ are disjoint intervals and:

\begin{enumerate}
\item the restriction of $T$ to each interval $I_{i}$ can be extended to a
$C^{1}$ monotonic map defined on $\overline{I_{i}}$ and the function
$h:I\rightarrow \mathbb{R}$ defined by $h(x)=|DT(x)|^{-1}$ has bounded
variation.

\item there are constants $C>0$ and $\sigma >1$ such that $\abs{DT^{n}(x)}>C\sigma ^{n}$ for every $n\geq 1$ and every $x\in I$ for which
the derivative is defined.

\item For each interval $J\subset I$ there is $n\geq 1$ such that $f^{n}(J)=I.$
\end{enumerate}

We remark that by point 1), in each interval $I_{i}$ the map is Lipschitz. We remark that  this restriction is not strictly necessary for what follows (see \cite{GalHoyRoj3}), we suppose it for the seek of simplicity.
As said before, by classical results this kind of map has an absolutely continuous invariant measure (see \cite{Via97}, chapter 3 e.g.).

\begin{theorem}
If $T$ a piecewise expanding map as above, then it has a unique ergodic
absolutely continuous invariant measure $\mu $. The basin of this SRB
measure has full Lebesgue measure. Moreover $\mu $ can be written as $\mu
=\phi m$ where $\phi $ has bounded variation and $m$ is the Lebesgue measure.
\end{theorem}

Moreover as before, the SRB measure is also computable

\begin{proposition}
If $T$ is a $m-$computable piecewise expanding map satisfying points
1),...,3) above then its SRB measure is computable.
\end{proposition}

\begin{proof}
Let us consider $\psi \in \mathcal{F}$ and $\integral{\psi\circ T^{n}}{m}$.
Since $T$ is $\delta$a.e computable for the Lesbegue measure, the ends of the intervals $I_{1},...,I_{k}$
are computable, and so is (uniformly on $n$) for $T^{n}$. If $T$ is
$l$-Lipschitz in each interval $I_{i}$ then $T^{n}$ is $l^{n}$-Lipschitz in
each of its continuity intervals and the Lipschitz constant of $\psi \circ
T^{n}$ on each continuity interval can be estimated. By this is
straightforward to show that also $\integral{\psi \circ T^{n}}{m}$ can be
calculated up to any accuracy (approximate $\psi \circ T^{n}$ by piecewise
constant functions and estimate error by Lipschitz constant).

Now, from \cite{Via97} proposition 3.8, remark 3.2 it holds that there are
$\lambda<1, C>0$ such that for each $\psi \in L^{1}$
\begin{equation*}
\abs{\integral{\psi \circ T^{n}}{m}-\integral{\psi}{\mu}}\leq C~\lambda^{n}~\norm{\psi}{L^1}.
\end{equation*}

This means that since $\integral{\psi \circ T^{n}}{m}$ can be calculated up
to any given error then the integral $\integral{\psi}{\mu}$ with respect to
$\mu $ can be calculated up to any given accuracy and by lemma \ref{gacs}
then $\mu $ is computable.
\end{proof}

As Unif. Hyperbolic systems, also Piecewise Expanding maps can be shown to
have exponential decay of correlations on bounded variation observables (see 
\cite{Via97} Remark 3.2) and BV norm of functions in $\mathcal{F}$ can be
estimated. Hence as in the previous section we obtain

\begin{corollary}
In an $m$-computable piecewise expanding system equipped with its SRB
measure, the set of computable typical points is dense in $[0,1]$.
\end{corollary}


\subsection{Manneville-Pomeau type maps}

We say that a map $T:[0,1]\rightarrow \lbrack 0,1]$ is a
\emph{Manneville-Pomeau type map (MP map)} with exponent $s$ if it
satisfies the following conditions:

\begin{enumerate}
\item there is $c\in (0,1)$ such that, if $I_{0}=[0,c]$ and $I_{1}=(c,1]$,
then $T\big|_{(0,c)}$ and $T\big|_{(c,1)}$ extend to $C^{1}$
diffeomorphisms, which is $C^{2}$ for $x>0,$ $T(I_{0})=[0,1]$,
$T(I_{1})=(0,1]$ and $T(0)=0$;

\item there is $\lambda >1$ such that $T^{\prime}\ge \lambda$ on $I_1$,
whereas $T^{\prime}>1$ on $(0,c]$ and $T^{\prime}(0)=1$;

\item the map $T$ has the following behaviour when $x\rightarrow 0^{+}$ 
\begin{equation*}
T(x)=x+rx^{1+s}(1+u(x))
\end{equation*}
for some constant $r>0$ and $s>0$ and $u$ satisfies $u(0)=0$ and $u^{\prime
}(x)=O(x^{t-1})$ as $x\rightarrow 0^{+}$ for some $t>0$.
\end{enumerate}

In \cite{Iso05} (see also \cite{Go}) it is proved that for $0<s<1$ these
systems have a unique absolutely continuous invariant measure, whose
density $f$ is locally Lipschitz in a neighborhood of each $x>0$ (the
density diverges at $x=0$) the system has polynomial decay of correlations
for $(1-s)$-H\"{o}lder observables. Moreover we have that:

\begin{theorem}
If $T$ is a computable MP map then its absolutely continuous invariant
measure $\mu$ is computable.
\end{theorem}

\begin{proof}
Let $f$ be the density of $\mu $. $T$ is topologically conjugated to the
doubling map $x\rightarrow 2x~(\func{mod}~1)~$ hence for each small
interval $I$ there is $k>0$ such that $T^{k}(I)=[0,1]$. Since $f$ is
locally Lipschitz, there is a small interval $J$ on which $f>\delta
_{1}>0$. Let $n$ be such that $T^{n}(J)=[0,1]$. Let $I$ be some small
interval, then there exist $J^{\prime }\subset J$ such that
$T^{n}(J^{\prime })=I$. Since $T$ is $\lambda $-Lipschitz, we have
$m(J^{\prime })\geq \frac{m(I)}{\lambda ^{n}}$.  By this, $\mu (J^{\prime
})\geq \frac{\delta _{1}m(I)}{\lambda ^{n}}$ and by the invariance of $\mu
$, $\frac{\mu (I)}{m(I)}\geq \frac{\delta _{1}}{\lambda ^{n}}$ and then,
as $I$ is arbitrary, for each $x\in \lbrack 0,1]$ we have
$f(x)>\frac{\delta _{1}}{\lambda ^{n}}>0$. In particular, $\frac{1}{f}$
is $(1-s)$-H\"{o}lder. Now we use the fact that the system has polynomial
decay of correlations for $(1-s)$-H\"{o}lder observables. Let us consider
$\phi \in \mathcal{F}$ then we have that $\frac{1}{f}d\mu =dm$ and
$\integral{\frac{1}{f}}{\mu}=1$, hence, by the decay of correlation of this
kind of maps
\begin{equation*}
\abs{\integral{\phi \circ T^{n}}{m}-\integral{\phi}{\mu}}=\abs{\integral{\phi \circ
  T^{n}\frac{1}{f}}{\mu} -\integral{\phi}{\mu} \integral{\frac{1}{f}}{\mu}
}\leq C\norm{\phi}{1-s}\norm{\frac{1}{f}}{1-s}n^{s-1}.
\end{equation*}
The norm $\norm{\phi}{1-s}$ can be
estimated for functions in $\mathcal{F}$, and then, as in the previous
examples we have a way to calculate $\integral{\phi}{\mu}$ for each $\phi \in 
\mathcal{F}$ and again by lemma \ref{gacs}, $\mu $ is computable.
\end{proof}

\begin{corollary}
In a computable Manneville-Pomeau type system, the set of computable typical
points is dense in $[0,1]$.
\end{corollary}

\section*{Acknowledgments}
The authors wish to thank Professors Giuseppe Longo and Stefano Marmi for scientific and logistic support during the elaboration of the present
paper.

\bibliographystyle{alpha}
\bibliography{bibliography}{}

\end{document}